\documentclass[11pt,reqno]{amsart}
\usepackage{amssymb,latexsym}
\usepackage{cite}
\usepackage[height=190mm,width=130mm]{geometry}

\theoremstyle{plain}
\newtheorem{theorem}{Theorem}
\newtheorem{lemma}{Lemma}
\newtheorem{corollary}{Corollary}

\theoremstyle{definition}

\theoremstyle{remark}

\numberwithin{equation}{section}
\input{tcilatex}

\begin{document}
\title[Positive Integer Solutions of the Pell Equation $x^{2}-dy^{2}=N,$ $%
d\in \left\{ k^{2}\pm 4,\text{ }k^{2}\pm 1\right\} $ \\
and $N\in \left\{ \pm 1,\pm 4\right\} $]{Positive Integer Solutions of the
Pell Equation $x^{2}-dy^{2}=N,$ $d\in \left\{ k^{2}\pm 4,\text{ }k^{2}\pm
1\right\} $ and $N\in \left\{ \pm 1,\pm 4\right\} $}
\author{Ref\.{i}k Kesk\.{i}n$^{\text{}}$}
\address{Sakarya University\\
Mathematics Department\\
TR 54187 Sakarya\\
Turkey}
\email{rkeskin@sakarya.edu.tr}
\author{Merve G\"{u}ney$^{\text{}}$}
\address{Sakarya University\\
Mathematics Department\\
TR 54187 Sakarya\\
Turkey}
\email{mer\_ney\_mat@hotmail.com}

\begin{abstract}
Let $\ k$ be a natural number and $d=k^{2}\pm 4$ or $k^{2}\pm 1$. In this
paper, by using continued fraction expansion of $\sqrt{d},$ we find
fundamental solution of the equations $x^{2}-dy^{2}=\pm 1$ and we get all
positive integer solutions of the equations $x^{2}-dy^{2}=\pm 1$ in terms of
generalized Fibonacci and Lucas sequences. Moreover, we find all positive
integer solutions of the equations $x^{2}-dy^{2}=\pm 4$ in terms of
generalized Fibonacci and Lucas sequences. Although some of the results are
well known, we think our method is elementary and different from the others.
\end{abstract}

\maketitle



\section{Introduction}

Let $d$ be a positive integer that is not a perfect square. It is well known
that the Pell equation $x^{2}-dy^{2}=1$ has always positive integer
solutions. When $N\neq 1,$ the Pell equation $x^{2}-dy^{2}=N$ may not have
any positive integer solutions. It can be seen that the equations $%
x^{2}-3y^{2}=-1$ and $x^{2}-7y^{2}=-4$ have no positive integer solutions.
Whether or not there exists a positive integer solution to the equation $%
x^{2}-dy^{2}=-1$ depends on the period length of the continued fraction
expansion of $\sqrt{d}$ (See section $2$ for more detailed information).
When $\ k$ is a positive integer and $d\in \left\{ k^{2}\pm 4,\text{ }%
k^{2}\pm 1\right\} ,$ positive integer solutions of the equations $%
x^{2}-dy^{2}=\pm 4$ and $x^{2}-dy^{2}=\pm 1$ have been investigated by Jones
in \cite{JONES} and the method used in the proofs of the theorems are method
of descent of Fermat. The same or similar equations are investigated by some
other authors in \cite{MCDANIEL}, \cite{MELHAM}, \cite{ZHIWEI}, \cite{KESKIN}%
, and \cite{ISMAIL}. Especially, when a solution exists, all positive
integer solutions of the equations $x^{2}-dy^{2}=\pm 4$ and $%
x^{2}-dy^{2}=\pm 1$ are given in terms of the generalized Fibonacci and
Lucas sequences. In this paper, if a solution exists, we will use continued
fraction expansion of $\sqrt{d}$ in order to get all positive integer
solutions of the equations $x^{2}-dy^{2}=\pm 1$ when $d\in \left\{ k^{2}\pm
4,\text{ }k^{2}\pm 1\right\} .$ Moreover, we will find all positive integer
solutions of the equations $x^{2}-dy^{2}=\pm 4$ when $d\in \left\{ k^{2}\pm
4,\text{ }k^{2}\pm 1\right\} .$ Our proofs are elementary and we think our
method is new and different from the others.

Now we briefly mention the generalized Fibonacci and Lucas sequences $\left(
U_{n}\left( k,s\right) \right) $ and $\left( V_{n}\left( k,s\right) \right) $%
. Let $k$ and $s$ be two nonzero integers with $k^{2}+4s>0.$ Generalized
Fibonacci sequence is defined by%
\begin{equation*}
U_{0}\left( k,s\right) =0,U_{1}\left( k,s\right) =1\text{ and }U_{n+1}\left(
k,s\right) =kU_{n}\left( k,s\right) +sU_{n-1}\left( k,s\right)
\end{equation*}%
for $n\geqslant 1$ and generalized Lucas sequence is defined by%
\begin{equation*}
V_{0}\left( k,s\right) =2,V_{1}\left( k,s\right) =k\text{ and }V_{n+1}\left(
k,s\right) =kV_{n}\left( k,s\right) +sV_{n-1}\left( k,s\right)
\end{equation*}%
for $n\geqslant 1,$ respectively. For $k=s=1,$ the sequences $\left(
U_{n}\right) $ and $\left( V_{n}\right) $ are called Fibonacci and Lucas
sequences and they are denoted as $\left( F_{n}\right) $ and $\left(
L_{n}\right) ,$ respectively. For $k=2$ and $s=1,$ the sequences $\left(
U_{n}\right) $ and $\left( V_{n}\right) $ are called Pell and Pell-Lucas
sequences and they are denoted as $\left( P_{n}\right) $ and $\left(
Q_{n}\right) ,$ respectively. It is well known that%
\begin{equation*}
U_{n}\left( k,s\right) =\dfrac{\alpha ^{n}-\beta ^{n}}{\alpha -\beta }\text{
and }V_{n}\left( k,s\right) =\alpha ^{n}+\beta ^{n}
\end{equation*}%
where $\alpha =\left( k+\sqrt{k^{2}+4s}\right) /2$ and $\beta =\left( k-%
\sqrt{k^{2}+4s}\right) /2.$ The above identities are known as Binet's
formulae. Clearly, $\alpha +\beta =k,$ $\alpha -\beta =\sqrt{k^{2}+4s},$ and
$\alpha \beta =-s.$ Especially, if $\alpha =\left( k+\sqrt{k^{2}+4}\right)
/2 $ and $\beta =\left( k-\sqrt{k^{2}+4}\right) /2$ , then we get%
\begin{equation}
U_{n}\left( k,1\right) =\dfrac{\alpha ^{n}-\beta ^{n}}{\alpha -\beta }\text{
and }V_{n}\left( k,1\right) =\alpha ^{n}+\beta ^{n}.  \label{1.2}
\end{equation}%
If $\alpha =\left( k+\sqrt{k^{2}-4}\right) /2$ and $\beta =\left( k-\sqrt{%
k^{2}-4}\right) /2$ , then we get%
\begin{equation}
U_{n}\left( k,-1\right) =\dfrac{\alpha ^{n}-\beta ^{n}}{\alpha -\beta }\text{
and }V_{n}\left( k,-1\right) =\alpha ^{n}+\beta ^{n}.  \label{1.3}
\end{equation}%
Also, if $\alpha =\left( 1+\sqrt{5}\right) /2$ and $\beta =\left( 1-\sqrt{5}%
\right) /2$ , then we get%
\begin{equation}
F_{n}=\dfrac{\alpha ^{n}-\beta ^{n}}{\alpha -\beta }\text{ and }L_{n}=\alpha
^{n}+\beta ^{n}.  \label{1.4}
\end{equation}%
Moreover, if $k$ is even, then it can be easily seen that%
\begin{equation*}
U_{n}\left( k,\pm 1\right) \text{ is odd}\Leftrightarrow n\text{ is odd},
\end{equation*}%
\begin{equation*}
U_{n}\left( k,\pm 1\right) \text{ is even}\Leftrightarrow n\text{ is even,}
\end{equation*}%
\begin{equation}
V_{n}\left( k,\pm 1\right) \text{ is even for all }n\in
\mathbb{N}
.  \label{1.8}
\end{equation}%
If $k$ is odd, then%
\begin{equation*}
2\mid V_{n}(k,\pm 1)\Leftrightarrow 2\mid U_{n}(k,\pm 1)\Leftrightarrow
3\mid n.
\end{equation*}%
For more information about generalized Fibonacci and Lucas sequences, one
can consult \cite{ROBINO}, \cite{KALMAN}, \cite{RIBENBOIM}, \cite{MCDANIEL},
and \cite{MELHAM}.

\section{Preliminaries}

Let $d$ be a positive integer which is not a perfect square and $N$ be any
nonzero fixed integer. Then the equation $x^{2}-dy^{2}=N$ is known as Pell
equation. For $N=\pm 1$, the equations $x^{2}-dy^{2}=1$ and $x^{2}-dy^{2}=-1$
are known as classical Pell equation. If $a^{2}-db^{2}=N$, we say that $%
(a,b) $ is a solution to the Pell equation $x^{2}-dy^{2}=N$. We use the
notations $(a,b)$ and $a+b\sqrt{d}$ interchangeably to denote solutions of
the equation $x^{2}-dy^{2}=N.$ Also, if $a$ and $b$ are both positive, we
say that $a+b\sqrt{d}$ is positive solution to the equation $x^{2}-dy^{2}=N.$%
Continued fraction plays an important role in solutions of the Pell
equations $x^{2}-dy^{2}=1$ and $x^{2}-dy^{2}=-1.$ Let $d$ be a positive
integer that is not a perfect square. Then there is a continued fraction
expansion of $\sqrt{d}$ such that $\sqrt{d}=\left[ a_{0},\overline{%
a_{1},a_{2},...,a_{l-1},2a_{0}}\right] $ where $l$ is the period length and
the $a_{j}$'s are given by the recursion formula;%
\begin{equation*}
\alpha _{0}=\sqrt{d},\text{ }a_{k}=\left\lfloor \alpha _{k}\right\rfloor
\text{ and }\alpha _{k+1}=\frac{1}{\alpha _{k}-a_{k}},\text{ }k=0,1,2,3,...
\end{equation*}%
Recall that $a_{l}=2a_{0}\ $and $a_{l+k}=a_{k}\ $for $k\geq 1$. The $n^{th}$
convergent of $\sqrt{d\text{ }}$ for $n\geq 0$ is given by%
\begin{equation*}
\frac{p_{n}}{q_{n}}=\left[ a_{0},a_{1},...,a_{n}\right] =a_{0}+\frac{1}{%
a_{1}+\frac{1}{\ddots \frac{1}{a_{n-1}+\frac{1}{a_{n}}}}}.\
\end{equation*}%
Let $x_{1}+y_{1}\sqrt{d}$ be a positive solution to the equation $%
x^{2}-dy^{2}=N$.\ We say that $x_{1}+y_{1}\sqrt{d}$ is the fundamental
solution of the equation $x^{2}-dy^{2}=N,$ if $x_{2}+y_{2}\sqrt{d}$ is
different solution to the equation$\ x^{2}-dy^{2}=N$, then $x_{1}+y_{1}\sqrt{%
d}<x_{2}+y_{2}\sqrt{d}$.Recall that if $a+b\sqrt{d}$ and $r+s\sqrt{d}$ are
two solutions to the equation $x^{2}-dy^{2}=N$, then $a=r$ if and only if $%
b=s,$ and $a+b\sqrt{d}<r+s\sqrt{d}$ if and only if $a<r$ and $b<s$. he
following lemmas and theorems can be found many elementary textbooks.

\begin{lemma}
\label{L1.1} If$\ x_{1}+y_{1}\sqrt{d}$ is the fundamental solution to the
equation $x^{2}-dy^{2}=-1$, then $(x_{1}+y_{1}\sqrt{d})^{2}$ is the
fundamental solution to the equation $x^{2}-dy^{2}=1$.\qquad \qquad \qquad
\qquad \qquad \qquad \qquad \qquad \qquad \qquad \qquad \qquad \qquad \qquad
\end{lemma}

If we know fundamental solution of the equations $x^{2}-dy^{2}=\pm 1$ and $%
x^{2}-dy^{2}=\pm 4$, then we can give all positive integer solutions to
these equations. Our theorems are as follows. For more information about
Pell equation, one can consult \cite{NAGELL}, \cite{JUDSON}, and \cite%
{ROBERTSON1}. Now we give the fundamental solution of the equations $%
x^{2}-dy^{2}=\pm 1$ by means of the period length of the continued fraction
expansion of $\sqrt{d}$.

\begin{lemma}
\label{L1.2} Let $l$ be the period length of continued fraction expansion of
$\sqrt{d}$. If $l$ is even, then the fundamental solution to the equation $%
x^{2}-dy^{2}=1$ is given by%
\begin{equation*}
x_{1}\bigskip +y_{1}\sqrt{d}=p_{l-1}+q_{l-1}\sqrt{d}
\end{equation*}%
and the equation $x^{2}-dy^{2}=-1$ has no integer solutions. If $l$ is odd,
then the fundamental solution of the equation $x^{2}-dy^{2}=1$ is given by$\
\ \ $%
\begin{equation*}
x_{1}\bigskip +y_{1}\sqrt{d}=p_{2l-1}+q_{2l-1}\sqrt{d}.
\end{equation*}%
and the fundamental solution to the equation $x^{2}-dy^{2}=-1$ is given by$\
\ \ $%
\begin{equation*}
x_{1}\bigskip +y_{1}\sqrt{d}=p_{l-1}+q_{l-1}\sqrt{d}.
\end{equation*}
\end{lemma}

\begin{theorem}
\label{t1.1} Let $x_{1}$\bigskip $+y_{1}\sqrt{d}$ be the fundamental
solution to the equation $x^{2}-dy^{2}=1$. Then all positive integer
solutions to the equation $x^{2}-dy^{2}=1$ are given by%
\begin{equation*}
x_{n}+y_{n}\sqrt{d}=(x_{1}+y_{1}\sqrt{d})^{n}
\end{equation*}%
with $n\geq 1$.
\end{theorem}

\begin{theorem}
\bigskip\ \label{t1.2} Let $x_{1}$\bigskip $+y_{1}\sqrt{d}$ be the
fundamental solution to the equation $x^{2}-dy^{2}=-1$. Then all positive
integer solutions to the equation $x^{2}-dy^{2}=-1$ are given by%
\begin{equation*}
x_{n}+y_{n}\sqrt{d}=(x_{1}+y_{1}\sqrt{d})^{2n-1}
\end{equation*}%
with $n\geq 1$. Now we give the following two theorems from \cite{JUDSON}.
See also \cite{ROBERTSON1}.
\end{theorem}

\begin{theorem}
\label{t1.3} Let $x_{1}$\bigskip $+y_{1}\sqrt{d}$ be the fundamental
solution to the equation $x^{2}-dy^{2}=4$. Then all positive integer
solutions to the equation $x^{2}-dy^{2}=4$ are given by%
\begin{equation*}
x_{n}+y_{n}\sqrt{d}=\frac{(x_{1}+y_{1}\sqrt{d})^{n}}{2^{n-1}}
\end{equation*}%
with $n\geq 1$.
\end{theorem}

\begin{theorem}
\label{t1.4} Let $x_{1}$\bigskip $+y_{1}\sqrt{d}$ be the fundamental
solution to the equation $x^{2}-dy^{2}=-4$. Then all positive integer
solutions to the equation $x^{2}-dy^{2}=-4$ are given by%
\begin{equation*}
x_{n}+y_{n}\sqrt{d}=\frac{(x_{1}+y_{1}\sqrt{d})^{2n-1}}{4^{n-1}}
\end{equation*}%
with $n\geq 1$.
\end{theorem}

From now on, we will assume that $k$ is a natural number. We give continued
fraction expansion of $\sqrt{d}$ for $d=k^{2}\pm 4.$ The proofs of the
following two theorems are easy and they can be found many text books on
number theory as an exercise (see, for example \cite{DON}).

\begin{theorem}
\label{t1.5} Let $k>1$. Then%
\begin{equation*}
\sqrt{k^{2}+4}=\left\{
\begin{array}{c}
\left[ k,\overline{\frac{k}{2},2k}\right] \text{ if }k\text{ is even,} \\
\left[ k,\overline{\frac{k-1}{2},1,1,\frac{k-1}{2},2k}\right] \text{ if }k%
\text{ is odd.}%
\end{array}%
\right.
\end{equation*}
\end{theorem}

\begin{theorem}
\label{t1.6} Let $k>3$. Then%
\begin{equation*}
\sqrt{k^{2}-4}=\left\{
\begin{array}{c}
\left[ k-1,\overline{1,\frac{k-3}{2},2,\frac{k-3}{2},1,2(k-1)}\right] \text{
if }k\text{ is odd,} \\
\left[ k-1,\overline{1,\frac{k-4}{2},1,2(k-1)}\right] \text{ if }k\text{ is
even and }k\neq 4 \\
\lbrack 3,\overline{2,6}]\text{ \ if }k=4%
\end{array}%
\right.
\end{equation*}
\end{theorem}

\begin{corollary}
\label{c1.1} Let $k>1$ and $d=k^{2}+4$. If $k$ is odd, then the fundamental
solution to the equation $x^{2}-dy^{2}=-1$ is
\begin{equation*}
x_{1}+y_{1}\sqrt{d}=\frac{k^{3}+3k}{2}+\frac{k^{2}+1}{2}\sqrt{d}.
\end{equation*}%
If $k$ is even, the equation $x^{2}-dy^{2}=-1$ has no positive integer
solutions.
\end{corollary}

\proof%
Assume that $k$ is odd. Then the period length of the continued fraction
expansion of $\sqrt{k^{2}+4}$ is $5$ by Theorem \ref{t1.5}. Therefore the
fundamental solution of the equation $x^{2}-dy^{2}=-1$ is $p_{4}+q_{4}\sqrt{d%
}$ by Lemma \ref{L1.2}. Since
\begin{equation*}
\frac{p_{4}}{q_{4}}=k+\frac{1}{\left( k-1\right) /2+\frac{1}{1+\frac{1}{1+%
\frac{1}{(k-1)/2}}}}=\frac{\frac{k^{3}+3k}{2}}{\frac{k^{2}+1}{2}},
\end{equation*}%
the proof follows. If $k$ is even, then the period length is even by Theorem %
\ref{t1.5} and therefore $x^{2}-dy^{2}=-1$ has no positive integer solutions
by Lemma \ref{L1.2}.%
\endproof%

\begin{corollary}
\label{c1.2} Let $k>1$ and $d=k^{2}+4$. Then the fundamental solution to the
equation $x^{2}-dy^{2}=1$ is
\begin{equation*}
x_{1}+y_{1}\sqrt{d}=\left\{
\begin{array}{c}
\frac{k^{2}+2}{2}+\frac{k}{2}\sqrt{d}\text{ \ if }k\text{ is even,} \\
\left( \frac{k^{3}+3k}{2}+\frac{k^{2}+1}{2}\sqrt[.]{d}\right) ^{2}\text{ }\
\text{if }k\text{ is odd}.%
\end{array}%
\right.
\end{equation*}
\end{corollary}

\proof%
If $k$ is even, then the proof follows from Lemma \ref{L1.2} and Theorem \ref%
{t1.5}. If $k$ is odd, then the proof follows from Corollary \ref{c1.1} and
Lemma \ref{L1.1}.%
\endproof%

From Lemma \ref{L1.2} and Theorem \ref{t1.6}, we can give the following
corollary.

\begin{corollary}
\label{c1.10} Let $k>3$ and $d=k^{2}-4$. Then the fundamental solution to
the equation $x^{2}-dy^{2}=1$ is given by%
\begin{equation*}
x_{1}+y_{1}\sqrt{d}=\left\{
\begin{array}{c}
\frac{k^{2}-2}{2}+\frac{k}{2}\sqrt{d}\text{ if }k\text{ is even,} \\
\frac{k^{3}-3k}{2}+\frac{k^{2}-1}{2}\sqrt{d}\text{ if }k\text{ is odd.}%
\end{array}%
\right. \text{ }
\end{equation*}
\end{corollary}

\begin{corollary}
\label{c1.7} Let $k>3$. Then the equation $x^{2}-(k^{2}-4)y^{2}=-1$ has no
integer solutions$.$
\end{corollary}

\proof%
The period length of continued fraction expansion of $\sqrt{k^{2}-4}$ is
always even by Theorem \ref{t1.6}. Thus by Lemma \ref{L1.2}, it follows that
there is no positive integer solutions of the equation $%
x^{2}-(k^{2}-4)y^{2}=-1$.%
\endproof%

\section{Main Theorems}

\begin{theorem}
\label{t1.19} Let $k>1$ and $d=k^{2}+4$. Then all positive integer solutions
of the equation $x^{2}-dy^{2}=1$ are given by%
\begin{equation*}
(x,y)=\left\{
\begin{array}{c}
\left( \frac{V_{2n}(k,1)}{2},\frac{U_{2n}(k,1)}{2}\right) \text{ if }k\text{
is even,} \\
\left( \frac{V_{6n}(k,1)}{2},\frac{U_{6n}(k,1)}{2}\right) \text{ if }k\text{
is odd},\
\end{array}%
\right.
\end{equation*}%
with $n\geq 1.$
\end{theorem}

\proof%
Assume that $k$ is even. Then, by Corollary \ref{c1.2} and Theorem \ref{t1.1}%
, all positive integer solutions of the equation $x^{2}-dy^{2}=1$ are given
by%
\begin{equation*}
x_{n}+y_{n}\sqrt{d}=\left( \frac{k^{2}+2}{2}+\frac{k}{2}\sqrt{d}\right) ^{n}
\end{equation*}%
with $n\geq 1$. Let $\alpha _{1}=\frac{k^{2}+2}{2}+\frac{k}{2}\sqrt{d}$ and $%
\beta _{1}=\frac{k^{2}+2}{2}-\frac{k}{2}\sqrt{d}$. Then
\begin{equation*}
x_{n}+y_{n}\sqrt{d}=\alpha _{1}^{n}\text{ and }x_{n}-y_{n}\sqrt{d}=\beta
_{1}^{n}.
\end{equation*}%
Thus it follows that $x_{n}=\frac{\alpha _{1}^{n}+\beta _{1}^{n}}{2}$ and $%
y_{n}=\frac{\alpha _{1}^{n}-\beta _{1}^{n}}{2\sqrt{d}}$. Let
\begin{equation*}
\alpha =\frac{k+\sqrt{k^{2}+4}}{2}\text{ \ and }\beta =\frac{k-\sqrt{k^{2}+4}%
}{2}.
\end{equation*}%
Then it is seen that $\alpha ^{2}=\alpha _{1}$ and $\beta ^{2}=\beta _{1}$.
Thus it follows that%
\begin{equation*}
x_{n}=\frac{\alpha ^{2n}+\beta ^{2n}}{2}=\frac{V_{2n}(k,1)}{2}
\end{equation*}%
and
\begin{equation*}
y_{n}=\frac{\alpha ^{2n}-\beta ^{2n}}{2\sqrt{d}}=\frac{\alpha ^{2n}-\beta
^{2n}}{2(\alpha -\beta )}=\frac{U_{2n}(k,1)}{2}
\end{equation*}%
by (\ref{1.2}). Now assume that \ $k$ is odd. Then by Corollary \ref{c1.2}
and Theorem \ref{t1.1}, we get
\begin{equation*}
x_{n}+y_{n}\sqrt{d}=\left( \left( \frac{k^{3}+3k}{2}+\frac{k^{2}+1}{2}\sqrt[.%
]{d}\right) ^{2}\right) ^{n}
\end{equation*}%
with $n\geq 1$. Let%
\begin{equation*}
\alpha _{1}=\left( \frac{k^{3}+3k}{2}+\frac{k^{2}+1}{2}\sqrt[.]{d}\right)
^{2}
\end{equation*}%
and%
\begin{equation*}
\beta _{1}=\left( \frac{k^{3}+3k}{2}-\frac{k^{2}+1}{2}\sqrt[.]{d}\right)
^{2}.
\end{equation*}%
Then
\begin{equation*}
x_{n}+y_{n}\sqrt{d}=\alpha _{1}^{n}\text{ and }x_{n}-y_{n}\sqrt{d}=\beta
_{1}^{n}.
\end{equation*}%
Thus it is seen that $x_{n}=\frac{\alpha _{1}^{n}+\beta _{1}^{n}}{2}$ and $%
y_{n}=\frac{\alpha _{1}^{n}-\beta _{1}^{n}}{2\sqrt{d}}$. Let
\begin{equation*}
\alpha =\frac{k+\sqrt{k^{2}+4}}{2}\text{ and }\beta =\frac{k-\sqrt{k^{2}+4}}{%
2}.
\end{equation*}%
Since $\alpha _{1}=\left( \frac{k^{3}+3k}{2}+\frac{k^{2}+1}{2}\sqrt[.]{d}%
\right) ^{2}=(\alpha ^{3})^{2}=\alpha ^{6}$ and thus $\beta _{1}=\beta ^{6}$%
, we get%
\begin{equation*}
x_{n}=\frac{\alpha ^{6n}+\beta ^{6n}}{2}=\frac{V_{6n}(k,1)}{2}
\end{equation*}%
and
\begin{equation*}
y_{n}=\frac{\alpha ^{6n}-\beta ^{6n}}{2\sqrt{d}}=\frac{\alpha ^{6n}-\beta
^{6n}}{2(\alpha -\beta )}=\frac{U_{6n}(k,1)}{2}
\end{equation*}%
by (\ref{1.2}). Then the proof follows.%
\endproof%

\begin{theorem}
\label{t1.21} Let $k>1$ be an odd integer and $d=k^{2}+4$. Then all positive
integer solutions of the equation $x^{2}-dy^{2}=-1$ are given by
\end{theorem}

\begin{equation*}
(x,y)=\left( \frac{V_{6n-3}(k,1)}{2},\frac{U_{6n-3}(k,1)}{2}\right) \ \
\end{equation*}%
with$\ n\geq 1.$

\proof%
Assume that $k>1$ be an odd integer. Then by Corollary \ref{c1.1} and
Theorem \ref{t1.2}, all positive integer solutions of the equation $%
x^{2}-dy^{2}=-1$ are given by
\begin{equation*}
x_{n}+y_{n}\sqrt{d}=\left( \frac{k^{3}+3k}{2}+\frac{k^{2}+1}{2}\sqrt{d}%
\right) ^{2n-1}
\end{equation*}%
with $n\geq 1$. Let $\alpha _{1}=\frac{k^{3}+3k}{2}+\frac{k^{2}+1}{2}\sqrt{d}
$ and $\beta _{1}=\frac{k^{3}+3k}{2}-\frac{k^{2}+1}{2}\sqrt{d}$. Then it
follows that

\begin{equation*}
x_{n}+y_{n}\sqrt{d}=\alpha _{1}^{2n-1}\text{ and }x_{n}-y_{n}\sqrt{d}=\beta
_{1}^{2n-1}
\end{equation*}%
and therefore $x_{n}=\frac{\alpha _{1}^{2n-1}+\beta _{1}^{2n-1}}{2}$ and $%
y_{n}=\frac{\alpha _{1}^{2n-1}-\beta _{1}^{2n-1}}{2\sqrt{d}}$. Let
\begin{equation*}
\alpha =\frac{k+\sqrt{k^{2}+4}}{2}\text{ and }\beta =\frac{k-\sqrt{k^{2}+4}}{%
2}.
\end{equation*}%
Then it is seen that
\begin{equation*}
\alpha ^{3}=\left( \frac{k+\sqrt{k^{2}+4}}{2}\right) ^{3}=\frac{k^{3}+3k}{2}+%
\frac{k^{2}+1}{2}\sqrt{d}=\alpha _{1}
\end{equation*}%
and
\begin{equation*}
\beta ^{3}=\left( \frac{k-\sqrt{k^{2}+4}}{2}\right) ^{3}=\frac{k^{3}+3k}{2}-%
\frac{k^{2}+1}{2}\sqrt{d}=\beta _{1}.
\end{equation*}%
Thus it follows that

\begin{equation*}
x_{n}=\frac{(\alpha ^{3})^{2n-1}+(\beta ^{3})^{2n-1}}{2}=\frac{\alpha
^{6n-3}+\beta ^{6n-3}}{2}=\frac{V_{6n-3}(k,1)}{2}
\end{equation*}%
and

\begin{equation*}
y_{n}=\frac{(\alpha ^{3})^{2n-1}-(\beta ^{3})^{2n-1}}{2\sqrt{d}}=\frac{%
\alpha ^{6n-3}-\beta ^{6n-3}}{2(\alpha -\beta )}=\frac{U_{6n-3}(k,1)}{2}
\end{equation*}%
by (\ref{1.2}). Then the proof follows.%
\endproof%

\begin{theorem}
\label{t1.22} Let $k>3$ and $d=k^{2}-4$. Then all positive integer solutions
of the equation $x^{2}-dy^{2}=1$ are given by%
\begin{equation*}
(x,y)=\left\{
\begin{array}{c}
\left( \frac{V_{2n}(k,-1)}{2},\frac{U_{2n}(k,-1)}{2}\right) \text{ if }k%
\text{ is even,} \\
\left( \frac{V_{3n}(k,-1)}{2},\frac{U_{3n}(k,-1)}{2}\right) \text{ if }k%
\text{ is odd,}%
\end{array}%
\right. \text{ }
\end{equation*}%
with $n\geq 1.$
\end{theorem}

\proof%
Assume that $k$ is even. By Corollary \ref{c1.10}\ and Theorem \ref{t1.1},
all positive integer solutions of the equation $x^{2}-dy^{2}=1$ are given by
\begin{equation*}
x_{n}+y_{n}\sqrt{d}=\left( \frac{k^{2}-2}{2}+\frac{k}{2}\sqrt{d}\right) ^{n}.
\end{equation*}%
Let $\alpha _{1}=\frac{k^{2}-2}{2}+\frac{k}{2}\sqrt{d}$ and $\beta _{1}=%
\frac{k^{2}-2}{2}-\frac{k}{2}\sqrt{d}.$ Then it follows that%
\begin{equation*}
x_{n}+y_{n}\sqrt{d}=\alpha _{1}^{n}\text{ and }x_{n}-y_{n}\sqrt{d}=\beta
_{1}^{n}.
\end{equation*}%
and therefore $x_{n}=\frac{\alpha _{1}^{n}+\beta _{1}^{n}}{2}$ and $y_{n}=%
\frac{\alpha _{1}^{n}-\beta _{1}^{n}}{2\sqrt{d}}$. Let
\begin{equation*}
\alpha =\frac{k+\sqrt{k^{2}-4}}{2}\text{ and }\beta =\frac{k-\sqrt{k^{2}-4}}{%
2}.
\end{equation*}%
Then it is seen that $\alpha ^{2}=\alpha _{1}$ and $\beta ^{2}=\beta _{1}$.
Thus it follows that%
\begin{equation*}
x_{n}=\frac{\alpha ^{2n}+\beta ^{2n}}{2}=\frac{V_{2n}(k,-1)}{2}\text{ }
\end{equation*}%
and%
\begin{equation*}
y_{n}=\frac{\alpha ^{2n}-\beta ^{2n}}{2\sqrt{d}}=\frac{\alpha ^{2n}-\beta
^{2n}}{2(\alpha -\beta )}=\frac{U_{2n}(k,-1)}{2}
\end{equation*}%
by (\ref{1.3})$.$ Now assume that $k$ is odd. Then by Corollary \ref{c1.10}
and Theorem \ref{t1.1}, we get

\begin{equation*}
x_{n}+y_{n}\sqrt{d}=\left( \frac{k^{3}-3k}{2}+\frac{k^{2}-1}{2}\sqrt{d}%
\right) ^{n}.
\end{equation*}%
Let $\alpha _{1}=\frac{k^{3}-3k}{2}+\frac{k^{2}-1}{2}\sqrt{d}$ and $\beta
_{1}=\frac{k^{3}-3k}{2}-\frac{k^{2}-1}{2}\sqrt{d}$. Then $x_{n}+y_{n}\sqrt{d}%
=\alpha _{1}^{n}$ and $x_{n}-y_{n}\sqrt{d}=\beta _{1}^{n}$. Thus it follows
that $x_{n}=\frac{\alpha _{1}^{n}+\beta _{1}^{n}}{2}$ and $y_{n}=\frac{%
\alpha _{1}^{n}-\beta _{1}^{n}}{2\sqrt{d}}$. Let $\alpha =\frac{k+\sqrt{%
k^{2}-4}}{2}$ and $\beta =\frac{k-\sqrt{k^{2}-4}}{2}.$ Since
\begin{equation*}
\alpha ^{3}=\left( \frac{k+\sqrt{k^{2}-4}}{2}\right) ^{3}=\frac{k^{3}-3k}{2}+%
\frac{k^{2}-1}{2}\sqrt{d}=\alpha _{1}
\end{equation*}%
and%
\begin{equation*}
\beta ^{3}=\left( \frac{k-\sqrt{k^{2}-4}}{2}\right) ^{3}=\frac{k^{3}-3k}{2}-%
\frac{k^{2}-1}{2}\sqrt{d}=\beta _{1},
\end{equation*}%
we get%
\begin{equation*}
x_{n}=\frac{\alpha ^{3n}+\beta ^{3n}}{2}=\frac{V_{3n}(k,-1)}{2}
\end{equation*}%
and%
\begin{equation*}
y_{n}=\frac{\alpha ^{3n}-\beta ^{3n}}{2\sqrt{d}}=\frac{\alpha ^{3n}-\beta
^{3n}}{2(\alpha -\beta )}=\frac{U_{3n}(k,-1)}{2}
\end{equation*}%
by (\ref{1.3}).%
\endproof%

Now we give all positive integer solutions of the equations $%
x^{2}-(k^{2}+4)y^{2}=\pm 4$ and $x^{2}-(k^{2}-4)y^{2}=\pm 4$. Before giving
all solutions of the equations $x^{2}-(k^{2}+4)y^{2}=\pm 4,$ we give the
following lemma which will be useful for finding the solutions.

\begin{lemma}
\label{L1.4} Let $a+b\sqrt{d}$ be a positive integer solution to equation $%
x^{2}-dy^{2}=4$. If $a>b^{2}-2$ , then $a+b\sqrt{d}$ is the fundamental
solution to the equation $x^{2}-dy^{2}=4$.
\end{lemma}

\proof%
If $b=1$, then the proof is trivial. Assume that $b>1$. Suppose that $%
x_{1}+y_{1}\sqrt{d}$ is a positive solution to the equation $x^{2}-dy^{2}=4$
such that $1\leq y_{1}<b$. Then it follows that $%
a^{2}-db^{2}=4=x_{1}^{2}-dy_{1}^{2}$ and thus $%
d=(x_{1}^{2}-4)/y_{1}^{2}=(a^{2}-4)/b^{2}$. This shows that $%
x_{1}^{2}b^{2}-y_{1}^{2}a^{2}=4b^{2}-4y_{1}^{2}=4(b^{2}-y_{1}^{2})>0$. Thus
\begin{equation*}
\lbrack \left( x_{1}b+y_{1}a\right) /2][\left( x_{1}b-y_{1}a\right)
/2]=b^{2}-y_{1}^{2}>1.
\end{equation*}%
It can be seen that $x_{1}b+y_{1}a$ and $x_{1}b-y_{1}a$ are even integers.
Let $(x_{1}b+y_{1}a)/2=k_{1}$ and $(x_{1}b-y_{1}a)/2=k_{2}$. Then $%
k_{1}k_{2}=b^{2}-y_{1}^{2}$ and $a=(k_{1}-k_{2})/y_{1}$. Thus
\begin{equation*}
a=\frac{k_{1}-k_{2}}{y_{1}}\leq \frac{k_{1}k_{2}-1}{y_{1}}=\frac{%
b^{2}-y_{1}^{2}-1}{y_{1}}\leq b^{2}-y_{1}^{2}-1\leq b^{2}-2,
\end{equation*}%
which is a contradiction since $a>b^{2}-2$. Then the proof follows.%
\endproof%

\begin{theorem}
\label{t1.7} Let $k>1$. Then all positive integer solutions of the equation $%
x^{2}-(k^{2}+4)y^{2}=4$ are given by%
\begin{equation*}
(x,y)=(V_{2n}(k,1),U_{2n}(k,1))\text{ }
\end{equation*}%
\bigskip with $n\geq 1.$
\end{theorem}

\proof%
Let $a=k^{2}+2$ and $b=k$. Then $a+b\sqrt{k^{2}+4}$ is a positive integer
solution of the equation $x^{2}-(k^{2}+4)y^{2}=4.$ Since $%
a=k^{2}+2>k^{2}-2=b^{2}-2,$ it follows that $k^{2}+2+k\sqrt{k^{2}+4}$ is the
fundamental solution of the equation $x^{2}-(k^{2}+4)y^{2}=4,$ by Lemma \ref%
{L1.4}. Thus, by Theorem \ref{t1.3}, all positive integer solutions of the
equation $x^{2}-dy^{2}=4$ are given by%
\begin{equation*}
x_{n}+y_{n}\sqrt{d}=\frac{(k^{2}+2+k\sqrt{k^{2}+4})^{n}}{2^{n-1}}=2\left(
\frac{k^{2}+2+k\sqrt{k^{2}+4}}{2}\right) ^{n}.
\end{equation*}%
Let $\alpha _{1}=\frac{k^{2}+2+k\sqrt{k^{2}+4}}{2}$ and $\beta _{1}=\frac{%
k^{2}+2-k\sqrt{k^{2}+4}}{2}$. Then it is seen that%
\begin{equation*}
x_{n}+y_{n}\sqrt{d}=2\alpha _{1}^{n}\text{ and }x_{n}-y_{n}\sqrt{d}=2\beta
_{1}^{n}.
\end{equation*}%
Thus it follows that $x_{n}=\alpha _{1}^{n}+\beta _{1}^{n}$ and $y_{n}=\frac{%
\alpha _{1}^{n}-\beta _{1}^{n}}{\sqrt{d}}$. Let
\begin{equation*}
\alpha =\frac{k+\sqrt{k^{2}+4}}{2}\text{ and }\beta =\frac{k-\sqrt{k^{2}+4}}{%
2}.
\end{equation*}%
Then
\begin{equation*}
\alpha ^{2}=\left( \frac{k+\sqrt{k^{2}+4}}{2}\right) ^{2}=\frac{k^{2}+2+k%
\sqrt{k^{2}+4}}{2}=\alpha _{1}
\end{equation*}%
and
\begin{equation*}
\beta ^{2}=\left( \frac{k-\sqrt{k^{2}+4}}{2}\right) ^{2}=\beta _{1}.
\end{equation*}%
Therefore we get%
\begin{equation*}
x_{n}=\alpha ^{2n}+\beta ^{2n}=V_{2n}(k,1)\text{ and }y_{n}=\frac{\alpha
^{2n}-\beta ^{2n}}{\sqrt{d}}=\frac{\alpha ^{2n}-\beta ^{2n}}{\alpha -\beta }%
=U_{2n}(k,1)
\end{equation*}%
by (\ref{1.2}). Then the proof follows.%
\endproof%

\begin{theorem}
\label{t1.23} Let $k>1$. Then all positive integer solutions of the equation
$x^{2}-(k^{2}+4)y^{2}=-4$ are given by%
\begin{equation*}
(x,y)=(V_{2n-1}(k,1),U_{2n-1}(k,1))
\end{equation*}%
with $n\geq 1.$
\end{theorem}

\proof%
Since $k^{2}-(k^{2}+4)=-4$, it follows that $k+\sqrt{k^{2}+4}$ is the
fundamental solution of the equation $x^{2}-(k^{2}+4)y^{2}=-4$. Thus by
Theorem \ref{t1.4}, all positive integer solutions of the equation $%
x^{2}-dy^{2}=-4$ are given by
\begin{equation*}
x_{n}+y_{n}\sqrt{d}=\frac{(k+\sqrt{k^{2}+4})^{2n-1}}{4^{n-1}}=2\left( \frac{%
k+\sqrt{k^{2}+4}}{2}\right) ^{2n-1}.
\end{equation*}%
Let $\alpha =\frac{k+\sqrt{k^{2}+4}}{2}$ and $\beta =\frac{k-\sqrt{k^{2}+4}}{%
2}$. Then it follows that
\begin{equation*}
x_{n}+y_{n}\sqrt{d}=2\alpha ^{2n-1}\text{ and }x_{n}-y_{n}\sqrt{d}=2\beta
^{2n-1}.
\end{equation*}%
Therefore%
\begin{equation*}
x_{n}=\alpha ^{2n-1}+\beta ^{2n-1}=V_{2n-1}(k,1)
\end{equation*}%
and%
\begin{equation*}
y_{n}=\frac{\alpha ^{2n-1}-\beta ^{2n-1}}{\sqrt{d}}=\frac{\alpha
^{2n-1}-\beta ^{2n-1}}{\alpha -\beta }=U_{2n-1}(k,1)
\end{equation*}%
by (\ref{1.2}). Then the proof follows.%
\endproof%

\begin{theorem}
\label{t1.24} Let $k>3$. Then all positive integer solutions of the equation
$x^{2}-(k^{2}-4)y^{2}=4$ are given by%
\begin{equation*}
(x,y)=(V_{n}(k,-1),U_{n}(k,-1))
\end{equation*}%
with $n\geq 1.$
\end{theorem}

\proof%
Since $k^{2}-(k^{2}-4)=4$, it is seen that $k+\sqrt{k^{2}-4}$ is the
fundamental solution of the equation $x^{2}-(k^{2}-4)y^{2}=4$. Thus, by
Theorem \ref{t1.3}$,$ all positive integer solutions of the equation $%
x^{2}-dy^{2}=4$ are given by
\begin{equation*}
x_{n}+y_{n}\sqrt{d}=\frac{(k+\sqrt{k^{2}-4})^{n}}{2^{n-1}}=2\left( \frac{k+%
\sqrt{k^{2}-4}}{2}\right) ^{n}.
\end{equation*}%
Let $\alpha =\frac{k+\sqrt{k^{2}-4}}{2}$ and $\beta =\frac{k-\sqrt{k^{2}-4}}{%
2}$. Then it follows that $x_{n}+y_{n}\sqrt{d}=2\alpha ^{n}$ and $x_{n}-y_{n}%
\sqrt{d}=2\beta ^{n}$. Thus we get%
\begin{equation*}
x_{n}=\alpha ^{n}+\beta ^{n}=V_{n}(k,-1)\text{ and }y_{n}=\frac{\alpha
^{n}-\beta ^{n}}{\sqrt{d}}=\frac{\alpha ^{n}-\beta ^{n}}{\alpha -\beta }%
=U_{n}(k,-1)
\end{equation*}%
by (\ref{1.3}).%
\endproof%

The following theorem is given in \cite{ROBERTSON2}.

\begin{theorem}
\label{t1.13} Let $d$ be odd positive integer. If the equation $%
x^{2}-dy^{2}=-4$ has positive integer solution, then the equation $%
x^{2}-dy^{2}=-1$ has positive integer solutions.
\end{theorem}

Now we give the continued fraction expansions of $\sqrt{k^{2}+1\text{ }}$and
$\sqrt{k^{2}-1}$. Since the continued fraction expansions of them are given
in \cite{JEAN}, we omit their proofs.

\begin{theorem}
\label{t1.9} Let $k\geq 1$. Then $\sqrt{k^{2}+1}=\left[ k,\overline{2k}%
\right] $ and if $k>1$, then $\sqrt{k^{2}-1}=\left[ k-1,\overline{1,2(k-1)}%
\right] $.
\end{theorem}

The proofs of the following corollaries follow from Lemma \ref{L1.2} and
Theorem \ref{t1.9} and therefore we omit their proofs.

\begin{corollary}
\label{c1.12} Let $k\geq 1$ and $d=k^{2}+1$. Then the fundamental solution
of the equation $x^{2}-dy^{2}=1$ is
\begin{equation*}
x_{1}+y_{1}\sqrt{d}=2k^{2}+1+2k\sqrt{d}.
\end{equation*}
\end{corollary}

\begin{corollary}
\label{c1.13} Let $k\geq 1$ and $d=k^{2}+1$. Then the fundamental solution
of the equation $x^{2}-dy^{2}=-1$ is%
\begin{equation*}
x_{1}+y_{1}\sqrt{d}=k+\sqrt{d}.
\end{equation*}
\end{corollary}

\begin{corollary}
\label{c1.14} Let $k>1$ and $d=k^{2}-1$. Then the fundamental solution of
the equation $x^{2}-dy^{2}=1$ is%
\begin{equation*}
x_{1}+y_{1}\sqrt{d}=k+\sqrt{d}.
\end{equation*}
\end{corollary}

\begin{theorem}
\label{t1.10} Let $k\geq 1$. Then all positive integer solutions of the
equation $x^{2}-(k^{2}+1)y^{2}=1$ are given by%
\begin{equation*}
(x,y)=\left( \frac{V_{2n}(2k,1)}{2},U_{2n}(2k,1)\right) \
\end{equation*}%
with $n\geq 1.$
\end{theorem}

\proof%
By Corollary \ref{c1.12} and Lemma \ref{L1.2}, it follows that all positive
integer solutions of the equation $x^{2}-(k^{2}+1)y^{2}=1$ are given by%
\begin{equation*}
x_{n}+y_{n}\sqrt{k^{2}+1}=\left( 2k^{2}+1+2k\sqrt{k^{2}+1}\right)
^{n}=\left( 2k^{2}+1+k\sqrt{(2k)^{2}+4}\right) ^{n}.
\end{equation*}%
Let $\alpha =\frac{2k+\sqrt{(2k)^{2}+4}}{2}$ and $\beta =\frac{2k-\sqrt{%
(2k)^{2}+4}}{2}$. Then
\begin{equation*}
\alpha ^{2}=\left( \frac{2k+\sqrt{(2k)^{2}+4}}{2}\right) ^{2}=2k^{2}+1+k%
\sqrt{(2k)^{2}+4}\text{ }
\end{equation*}%
and

\begin{equation*}
\beta ^{2}=\left( \frac{2k-\sqrt{(2k)^{2}+4}}{2}\right) ^{2}=2k^{2}+1-k\sqrt{%
(2k)^{2}+4}.
\end{equation*}%
Thus it follows that
\begin{equation*}
x_{n}+y_{n}\sqrt{k^{2}+1}=x_{n}+\frac{y_{n}}{2}\sqrt{(2k)^{2}+4}=\alpha ^{2n}%
\text{ }
\end{equation*}%
and
\begin{equation*}
x_{n}-y_{n}\sqrt{k^{2}+1}=x_{n}-\frac{y_{n}}{2}\sqrt{(2k)^{2}+4}=\beta ^{2n}.
\end{equation*}%
Then it is seen that%
\begin{equation*}
x_{n}=\frac{\alpha ^{2n}+\beta ^{2n}}{2}=\frac{V_{2n}(2k,1)}{2}
\end{equation*}%
and%
\begin{equation*}
y_{n}=\frac{\alpha ^{2n}-\beta ^{2n}}{\sqrt{(2k)^{2}+4}}=\frac{\alpha
^{2n}-\beta ^{2n}}{\alpha -\beta }=U_{2n}(2k,1)
\end{equation*}%
by (\ref{1.2}).%
\endproof%

Since the proof of the following theorems are similar to that of above
theorems, we omit them.

\begin{theorem}
\label{t1.11} Let $k\geq 1$. Then all positive integer solutions of the
equation $x^{2}-(k^{2}+1)y^{2}=-1$ are given by%
\begin{equation*}
(x,y)=\left( \frac{V_{2n-1}(2k,1)}{2},U_{2n-1}(2k,1)\right)
\end{equation*}%
with $n\geq 1.$
\end{theorem}

\begin{theorem}
\label{t1.12} Let $k>1$. Then all positive integer solutions of the equation
$x^{2}-(k^{2}-1)y^{2}=1$ are given by%
\begin{equation*}
(x,y)=\left( \frac{V_{n}(2k,-1)}{2},U_{n}(2k,-1)\right)
\end{equation*}%
with $n\geq 1.$
\end{theorem}

\begin{corollary}
\label{c1.15} Let $k>1$. Then the equation $x^{2}-(k^{2}-1)y^{2}=-1$ has no
positive integer solutions.
\end{corollary}

\proof%
The period length of continued fraction expansion of $\sqrt{k^{2}-1}$ is
always even by Theorem \ref{t1.9}. Thus, by Lemma \ref{L1.2}, it follows
that there is no positive integer solutions of the equation $%
x^{2}-(k^{2}-1)y^{2}=-1$.%
\endproof%

\begin{theorem}
\bigskip \label{t1.8} Let $k>3$. Then the equation $x^{2}-(k^{2}-4)y^{2}=-4$
has no positive integer solutions.
\end{theorem}

\proof%
Assume that $k$ is odd. Then $k^{2}-4$ is odd and thus the proof follows
from Theorem \ref{t1.13} and Corollary \ref{c1.7}. Now assume that $k$ is
even. If $(a,b)$ is a solution to the equation $x^{2}-(k^{2}-4)y^{2}=-4,$
then $a$ is even. Thus we get%
\begin{equation*}
\left( a/2\right) ^{2}-((k/2)^{2}-1)b^{2}=-1,
\end{equation*}%
which is impossible by Corollary \ref{c1.15}. Then the proof follows.%
\endproof%

Now we give all positive integer solutions of the equations $%
x^{2}-(k^{2}+1)y^{2}=\pm 4$ and $x^{2}-(k^{2}-1)y^{2}=\pm 4$.

\begin{theorem}
\label{t1.18} Let $k\geq 1$ and $k\neq 2$. Then all positive integer
solutions of the equation $x^{2}-(k^{2}+1)y^{2}=-4\ $are given by
\begin{equation*}
(x,y)=\left( V_{2n-1}(2k,1),2U_{2n-1}(2k,1)\right)
\end{equation*}%
with $n\geq 1.$
\end{theorem}

\proof%
Since $k\geq 1$, it can be shown that $2k+2\sqrt{k^{2}+1}$ is the
fundamental solution to the equation $x^{2}-(k^{2}+1)y^{2}=-4.$ Then by
Theorem \ref{t1.4}, all positive integer solutions of the equation $%
x^{2}-(k^{2}+1)y^{2}=-4$ are given by
\begin{equation*}
x_{n}+y_{n}\sqrt{k^{2}+1}=2\left( \frac{2k+2\sqrt{k^{2}+1}}{2}\right)
^{2n-1}=2\left( \frac{2k+\sqrt{(2k)^{2}+4}}{2}\right) ^{2n-1}.
\end{equation*}%
Let $\alpha =\frac{2k+\sqrt{(2k)^{2}+4}}{2}$ and $\beta =\frac{2k-\sqrt{%
(2k)^{2}+4}}{2}.$ Then we get
\begin{equation*}
x_{n}+y_{n}\sqrt{k^{2}+1}=x_{n}+\frac{y_{n}}{2}\sqrt{(2k)^{2}+4}=2\alpha
^{2n-1}
\end{equation*}%
and%
\begin{equation*}
x_{n}-y_{n}\sqrt{k^{2}+1}=x_{n}-\frac{y_{n}}{2}\sqrt{(2k)^{2}+4}=2\beta
^{2n-1}.
\end{equation*}%
Thus it follows that
\begin{equation*}
x_{n}=\alpha ^{2n-1}+\beta ^{2n-1}=V_{2n-1}(2k,1)
\end{equation*}%
and
\begin{equation*}
y_{n}=2\frac{\alpha ^{2n-1}-\beta ^{2n-1}}{\sqrt{(2k)^{2}+4}}=2\frac{\alpha
^{2n-1}-\beta ^{2n-1}}{\alpha -\beta }=2U_{2n-1}(2k,1)
\end{equation*}%
by (\ref{1.2}).%
\endproof%

Now we can give the following corollary from Theorem\ \ref{t1.18} and
identity (\ref{1.8}).

\begin{corollary}
\label{c1.28} If $\left( a,b\right) $ is a positive integer solution of the
equation $x^{2}-(k^{2}+1)y^{2}=-4,$ then $a$ and $b$ are even.
\end{corollary}

Since the proof of the following theorem is similar to that of Theorem \ref%
{t1.18} , we omit it.

\begin{theorem}
\label{t1.17} Let $k>1$. Then all positive integer solutions of the equation
$x^{2}-(k^{2}-1)y^{2}=4$ are given by
\begin{equation*}
(x,y)=\left( V_{n}(2k,-1),2U_{n}(2k,-1)\right)
\end{equation*}%
with $n\geq 1.$
\end{theorem}

\begin{theorem}
\label{t1.20} Let $k\geq 1$ and $k\neq 2$. Then all positive integer
solutions of the equation $x^{2}-(k^{2}+1)y^{2}=4$ are given by%
\begin{equation*}
(x,y)=\left( V_{2n}(2k,1),2U_{2n}(2k,1)\right)
\end{equation*}%
with $n\geq 1.$
\end{theorem}

\proof%
Firstly, we show that if $(a,b)$ is a solution to the equation $%
x^{2}-(k^{2}+1)y^{2}=4$, then $a$ and $b$ are even. Assume that $k$ is odd.
Then $k^{2}+1=2t$ for some odd integer $t.$ Since $a^{2}-2tb^{2}=4$, it
follows that $a$ is even and therefore $b$ is even. Now assume that $k$ is
even. Let $d=k^{2}+1.$ Then $d$ is odd. Assume that $a$ and $b$ are odd
integers$.$ Let $x_{1}=\left\vert db-ka\right\vert ,$ $y_{1}=\left\vert
a-kb\right\vert .$ Then $x_{1}$ and $y_{1}$ are odd integers. Moreover,%
\begin{equation*}
x_{1}^{2}-dy_{1}^{2}=(db-ka)^{2}-d(a-kb)^{2}=b^{2}d(d-k^{2})+a^{2}(k^{2}-d)=b^{2}d-a^{2}\bigskip =-(a^{2}-db^{2})=-4
\end{equation*}%
Thus $x_{1}+y_{1}\sqrt{d}$ is a positive solution of the equation $%
x^{2}-(k^{2}+1)y^{2}=-4,$ which is impossible by Corollary \ref{c1.28}.
Therefore if $a+b\sqrt{d}$ is any solutions of the equation $x^{2}-dy^{2}=4,$
then $a$ and $b$ are even integers and thus $\frac{a}{2}+\frac{b}{2}\sqrt{d%
\text{ }}$ is a solution to the equation $x^{2}-dy^{2}=1$. Then it follows
that the fundamental solution of the equation $x^{2}-dy^{2}=4$ is $%
4k^{2}+2+4k\sqrt{d},$ by corollary \ref{c1.12}. Thus by Theorem \ref{t1.3},
it follows that all positive integer solutions of the equation $%
x^{2}-(k^{2}+1)y^{2}=4$ are given by
\begin{equation*}
x_{n}+y_{n}\sqrt{k^{2}+1}=2\left( \frac{4k^{2}+2+4k\sqrt{k^{2}+1}}{2}\right)
^{n}=2\left( \frac{4k^{2}+2+2k\sqrt{(2k)^{2}+4}}{2}\right) ^{n}.
\end{equation*}%
Let $\alpha =\frac{2k+\sqrt{(2k)^{2}+4}}{2}$ and $\beta =\frac{2k-\sqrt{%
(2k)^{2}+4}}{2}$. Then
\begin{equation*}
\alpha ^{2}=\left( \frac{2k+\sqrt{(2k)^{2}+4}}{2}\right) ^{2}=\frac{%
4k^{2}+2+2k\sqrt{(2k)^{2}+4}}{2}
\end{equation*}%
and

\begin{equation*}
\beta ^{2}=\left( \frac{2k-\sqrt{(2k)^{2}+4}}{2}\right) ^{2}=\frac{%
4k^{2}+2-2k\sqrt{(2k)^{2}+4}}{2}.
\end{equation*}%
Thus it follows that $x_{n}+y_{n}\sqrt{k^{2}+1}=x_{n}+\frac{y_{n}}{2}\sqrt{%
(2k)^{2}+4}=2\alpha ^{2n}$ and $x_{n}-\frac{y_{n}}{2}\sqrt{(2k)^{2}+4}%
=2\beta ^{2n}$. Then it is seen that
\begin{equation*}
x_{n}=\alpha ^{2n}+\beta ^{2n}=V_{2n}(2k,1)
\end{equation*}%
and
\begin{equation*}
y_{n}=2\frac{\alpha ^{2n}-\beta ^{2n}}{\sqrt{(2k)^{2}+4}}=2\frac{\alpha
^{2n}-\beta ^{2n}}{\alpha -\beta }=2U_{2n}(2k,1),
\end{equation*}%
by (\ref{1.2}).%
\endproof%

It can be shown that if $k>2,$ then the continued fraction expansion of $%
\sqrt{k^{2}-k}$ is $[k-1,\overline{1,2(k-1)}]$ (see \cite{DON}, page 234).
Therefore we can give the following corollary easily.

\begin{corollary}
\label{c1.22} Let $k>2$. Then the equation $x^{2}-(k^{2}-k)y^{2}=-1$ has no
positive integer solutions.
\end{corollary}

\begin{corollary}
\label{c1.9} Let $k\geq 2$ and $k\neq 3$. Then the equation $%
x^{2}-(k^{2}-1)y^{2}=-4$ has no positive integer solutions.
\end{corollary}

\proof%
Assume that $k$ is even. Then $k^{2}-1$ is odd and the proof follows from
Theorem \ref{t1.13} and Corollary \ref{c1.15}.

Assume that $k$ is odd and . Then $k^{2}-1$ is even. Now assume that $%
a^{2}-(k^{2}-1)b^{2}=-4$ for some positive integers $a$ and $b$. Then $a$ is
even and this implies that
\begin{equation*}
\left( a/2\right) ^{2}-[\left( k^{2}-1\right) /4]b^{2}=-1.
\end{equation*}%
This is impossible by Corollary \ref{c1.22}, since
\begin{equation*}
\left( k^{2}-1\right) /4=((k+1)/2)^{2}-(k+1)/2.
\end{equation*}%
\endproof%

Continued fraction expansion of $\sqrt{5}$ is $\left[ 2,\overline{4}\right]
. $ Then the period length of the continued fraction expasion of $\sqrt{5}$
is 1. Therefore the fundamental solution to the equation $x^{2}-5y^{2}=1$ is
$9+4\sqrt{5}$ and the fundamental solution to the equation $x^{2}-5y^{2}=-1$
is $2+\sqrt{5}$ by Lemma $\ref{L1.2}$. Therefore, by using ($\ref{1.4})$, we
can give the following corollaries easily.

\begin{corollary}
\label{c1.16} All positive integer solutions of the equation $x^{2}-5y^{2}=1$
are given by
\begin{equation*}
(x,y)=(\frac{L_{6n}}{2},\frac{F_{6n}}{2})
\end{equation*}%
~
\end{corollary}

with $n\geq 1.$

\begin{corollary}
\label{c1.17} All positive integer solutions of the equation $%
x^{2}-5y^{2}=-1 $ are given by
\begin{equation*}
(x,y)=(\frac{L_{6n-3}}{2},\frac{F_{6n-3}}{2})
\end{equation*}%
~
\end{corollary}

with $n\geq 1.$

It can be seen that fundamental solutions of the equations $x^{2}-5y^{2}=-4$
and $x^{2}-5y^{2}=4$ are $1+\sqrt{5}$ and $3+\sqrt{5},$ respectively. Thus
we can give following corollaries.

\begin{corollary}
\label{c1.20} All positive integer solutions of the equation $x^{2}-5y^{2}=4$
are given by%
\begin{equation*}
(x,y)=(L_{2n},F_{2n})
\end{equation*}%
~
\end{corollary}

with $n\geq 1.$

\begin{corollary}
\label{c1.21} All positive integer solutions of the equation $%
x^{2}-5y^{2}=-4 $ are given by
\begin{equation*}
(x,y)=(L_{2n-1},F_{2n-1})
\end{equation*}%
~
\end{corollary}

with $n\geq 1.$


\begin{thebibliography}{99}
\bibitem{ADLER} Adler, A. and Coury, J. E., \emph{The Theory of Numbers: A
Text and Source Book of Problems,} Jones and Bartlett Publishers, Boston,
MA, 1995.

\bibitem{DON} Don Redmond, \emph{Number Theory: An Introduction}, Markel
Dekker, Inc, 1996

\bibitem{JEAN} Jean-Marie De Koninck and Armel Mercier, \emph{1001 Problems
in Classical Number Theory}, American Mathematical Society, 2007.

\bibitem{ROBERTSON1} John P. Robertson, \emph{Solving the generalized Pell
equation }$x^{2}-Dy^{2}=N$. http://hometown.aol.com/jpr2718/pell.pdf, May
2003. (Description of LMM Algorithm for solving Pell's equation).

\bibitem{ROBERTSON2} John P. Robetson, \emph{On D so that }$x^{2}-Dy^{2}$%
\emph{\ represents }$m$\emph{\ and }$-m$\emph{\ and not }$-1,$ Acta
Mathematica Academia Paedogogocae Nyiregyhaziensis, 25 (2009) 155-164.

\bibitem{JONES} J. P. Jones, \emph{Representation of Solutions of Pell
Equations Using Lucas Sequences, }Acta Academia Pead. Agr., Sectio
Mathematicae 30 (2003) 75-86.

\bibitem{KALMAN} D. Kalman, R. Mena, \emph{The Fibonacci Numbers exposed},
Mathematics Magazine 76(2003) 167-181.

\bibitem{KESKIN} R. Keskin, \emph{Solutions of some quadratic Diophantine
equations}, Computers and Mathematics with Applications, 60 (2010) 2225-2230.

\bibitem{MCDANIEL} W. L. McDaniel, \emph{Diophantine Representation of Lucas
Sequences, }The Fibonacci Quarterly 33 (1995), 58-63.

\bibitem{MELHAM} R. Melham, \emph{Conics Which Characterize Certain Lucas
Sequences, }The Fibonacci Quarterly 35 (1997), 248-251.

\bibitem{WILLIAMS} Michael J. Jacobson, Hugh C. Williams, \emph{Solving the
Pell Equation}, Springer, (2006).

\bibitem{NAGELL} T. Nagell, \emph{Introduction to Number Theory, }Chelsea
Publishing Company, New York, 1981.

\bibitem{RIBENBOIM} P. Ribenboim, \emph{My Numbers, My Friends, }%
Springer-Verlag New York, Inc., 2000.

\bibitem{ROBINO} S. Robinowitz, \emph{Algorithmic Manipulation of Fibonacci
Identities, in: Application of Fibonacci Numbers}, vol. 6, Kluwer Academic
Pub., Dordrect, The Netherlands, 1996, pp. 389-408.

\bibitem{JUDSON} J. LeVeque, W., \emph{Topics in Number Theory,} Volume 1
and 2, Dover Publications 2002.

\bibitem{ISMAIL} M. E. H. Ismail, \emph{One Parameter Generalizations of the
Fibonacci and Lucas Numbers}, The Fibonacci Quarterly 46-47 (2009), 167-180.

\bibitem{ZHIWEI} Zhiwei, S., \emph{Singlefold Diophantine Representation of
the Sequence }$u_{0}=0,u_{1}=1$\emph{\ and }$u_{n+2}=mu_{n+1}+u_{n},$ Pure
and Applied Logic, Beijing Univ. Press, Beijing, 97-101, 1992.
\end{thebibliography}
\end{document}